\documentclass[12pt, leqno]{amsart}

\usepackage{texdraw,multicol,rotating}

\usepackage{amsmath}
\usepackage{amssymb}
\usepackage{enumerate}

\input{txdtools.tex}

\numberwithin{equation}{section}

\theoremstyle{plain}
\newtheorem{thm}{Theorem}[section]
\newtheorem{prop}[thm]{Proposition}
\newtheorem{lem}[thm]{Lemma}
\newtheorem{cor}[thm]{Corollary}

\theoremstyle{definition}

\newtheorem{example}[thm]{Example}
\newtheorem{remark}[thm]{Remark}

\theoremstyle{remark}

\newbox{\tmpa}
\newbox{\tmpb}

\DeclareMathOperator{\wt}{wt}
\newcommand{\nc}{\newcommand}
\nc{\Uq}{U_q} \nc{\Z}{\mathbf{Z}} \nc{\C}{\mathbf{C}}
\nc{\Q}{\mathbf{Q}}
\nc{\op}{\oplus} \nc{\ot}{\otimes} \nc{\pv}{P^{\vee}}
\nc{\ali}{\alpha_i} \nc{\B}{\mathbf{B}} \nc{\F}{\mathbf{F}}
\nc{\bP}{\mathbf{P}} \nc{\V}{\mathbf{V}} \nc{\La}{\Lambda}
\nc{\la}{\lambda} \nc{\nbinom}[2]{\genfrac{}{}{0pt}{1}{{#1}}{{#2}}}
\nc{\qbinom}[2]{\left[\genfrac{}{}{0pt}{1}{{#1}}{{#2}}\right]}
\nc{\path}{\mathcal{P}} \nc{\fit}{\tilde{f}_i}
\nc{\eit}{\tilde{e}_i} \nc{\fjt}{\tilde{f}_j} \nc{\ejt}{\tilde{e}_j}
\nc{\Y}{\mathbf{Y}} \nc{\A}{\mathbf{A}} \nc{\ra}{\rightarrow}
\nc{\vep}{\varepsilon} \nc{\vphi}{\varphi} \nc{\vp}{\varphi}
\nc{\g}{\mathfrak{g}} \nc{\h}{\mathfrak{h}} \nc{\oP}{\overline{P}}
\nc{\pathp}{\mathbf{p}}

\nc{\tris}{ \bsegment \move(0 0)\lvec(10 0)\lvec(10 10)\lvec(0
0)\ifill f:0.7 \esegment } \nc{\recs}{ \bsegment \move(0 0)\lvec(10
0)\lvec(10 5)\lvec(0 5)\lvec(0 0)\ifill f:0.7 \esegment }
\nc{\hcvec}[5]{%
\getpos(#1 #3)\spx\spy \getpos(#2 #3)\epx\epy \getpos(#4
#5)\xoff\yoff \realadd \spx \xoff \twox \realadd \epx {-\xoff} \thrx
\realadd \spy \yoff \posy \move({\spx} {\spy}) \clvec ({\twox}
{\posy})({\thrx} {\posy})({\epx} {\epy}) \rmove(0 0) }

\nc{\ahead}[2]{%
\cossin (0 0)({#1} {#2})\cosa\sina \bsegment
  \drawdim in \setunitscale 0.065
  \realmult {-0.5} \cosa \hcosa
  \realmult {-0.5} \sina \hsina
  \move({\hcosa} {\hsina}) \ravec({\cosa} {\sina})
\esegment }
\nc{\boxi}{%
{%
\savebox{\tmppic}{\begin{texdraw} \small \drawdim em \textref h:C
v:C \setunitscale 0.55 \htext(0 0){$i$} \move(-1 -1)\lvec(-1
1)\lvec(1 1)\lvec(1 -1)\lvec(-1 -1)
\end{texdraw}}%
\raisebox{-0.19\height}{\usebox{\tmppic}}%
}%
}
\nc{\boxj}{%
{%
\savebox{\tmppic}{\begin{texdraw} \small \drawdim em \textref h:C
v:C \setunitscale 0.55 \htext(0 0.1){$j$} \move(-1 -1)\lvec(-1
1)\lvec(1 1)\lvec(1 -1)\lvec(-1 -1)
\end{texdraw}}%
\raisebox{-0.19\height}{\usebox{\tmppic}}%
}%
}
\nc{\boxipo}{%
{%
\savebox{\tmppic}{\begin{texdraw} \small \drawdim em \textref h:C
v:C \setunitscale 0.55 \htext(0.15 0){$i\!\!+\!\!1$} \move(-1.4
-1)\lvec(-1.4 1)\lvec(1.4 1)\lvec(1.4 -1)\lvec(-1.4 -1)
\end{texdraw}}%
\raisebox{-0.19\height}{\usebox{\tmppic}}%
}%
} \everytexdraw{ \drawdim in \arrowheadsize l:0.065 w:0.03
\arrowheadtype t:F \textref h:C v:C }
\newsavebox{\tmppic}
\newsavebox{\tmpfig}
\newsavebox{\tmpdraw}
\newsavebox{\tmpfiga}
\newsavebox{\tmpfigb}
\newsavebox{\tmpfigc}
\newsavebox{\tmpfigd}
\newsavebox{\tmpfige}
\newsavebox{\tmpfigf}
\newsavebox{\tmpfigg}
\newsavebox{\tmpfigh}
\newsavebox{\tmpfigi}
\newsavebox{\tmpfigj}
\newsavebox{\tmpfigk}
\newsavebox{\tmpfigl}
\newsavebox{\tmpfigm}
\newsavebox{\tmpfign}
\newsavebox{\tmpfigo}
\newsavebox{\tmpfigp}
\newsavebox{\tmpfigq}
\newsavebox{\tmpfigr}
\newsavebox{\tmpfigs}
\newsavebox{\tmpfigt}
\newsavebox{\tmpfigu}
\newsavebox{\tmpfigv}
\newsavebox{\tmpfigw}
\newsavebox{\tmpfigx}
\newsavebox{\tmpfigy}
\newsavebox{\tmpfigz}

\newsavebox{\tmpfigaa}
\newsavebox{\tmpfigab}
\newsavebox{\tmpfigac}
\newsavebox{\tmpfigad}
\newsavebox{\tmpfigae}
\newsavebox{\tmpfigaf}
\newsavebox{\tmpfigag}
\newsavebox{\tmpfigah}
\newsavebox{\tmpfigai}
\newsavebox{\tmpfigaj}
\newsavebox{\tmpfigak}
\newsavebox{\tmpfigal}
\newsavebox{\tmpfigam}
\newsavebox{\tmpfigan}
\newsavebox{\tmpfigao}
\newsavebox{\tmpfigap}
\newsavebox{\tmpfigaq}
\newsavebox{\tmpfigar}
\newsavebox{\tmpfigas}
\newsavebox{\tmpfigat}
\newsavebox{\tmpfigau}
\newsavebox{\tmpfigav}
\newsavebox{\tmpfigaw}
\newsavebox{\tmpfigax}
\newsavebox{\tmpfigay}
\newsavebox{\tmpfigaz}

\newsavebox{\tmpfigba}
\newsavebox{\tmpfigbb}
\newsavebox{\tmpfigbc}
\newsavebox{\tmpfigbd}
\newsavebox{\tmpfigbe}
\newsavebox{\tmpfigbf}
\newsavebox{\tmpfigbg}
\newsavebox{\tmpfigbh}

\nc{\node}{\lcir r:1 }
\nc{\sline}{\bsegment\savepos(10 0)(*ex *ey)
            \move(1 0)\rlvec(8 0)
            \esegment\move(*ex *ey)}
\nc{\dline}{\bsegment\savepos(10 0)(*ex *ey)
            \move(0.93 0.4)\rlvec(8.14 0)\rmove(0 -0.8)\rlvec(-8.14 0)
            \esegment\move(*ex *ey)}
\nc{\uline}{\bsegment\savepos(0 10)(*ex *ey)
            \move(0 1)\rlvec(0 8)
            \esegment\move(*ex *ey)}
\nc{\lpoint}{\savecurrpos(*ex *ey)
             \rmove(2.5 1.5)\rlvec(-1.5 -1.5)\rlvec(1.5 -1.5)
             \move(*ex *ey)}
\nc{\rpoint}{\savecurrpos(*ex *ey)
             \rmove(-2.5 -1.5)\rlvec(1.5 1.5)\rlvec(-1.5 1.5)
             \move(*ex *ey)}
\nc{\bline}{\bsegment\savepos(10 0)(*ex *ey)
            \linewd 0.6 \move(1.1 0)\rlvec(7.8 0)
            \esegment\move(*ex *ey)}

\nc{\araise}[1]{\raisebox{4.5pt}{#1}}
\nc{\braise}[1]{\raisebox{12.1pt}{#1}}
\nc{\craise}[1]{\raisebox{8pt}{#1}}
\nc{\draise}[1]{\raisebox{12pt}{#1}} \nc{\be}{\begin{enumerate}}
\nc{\ee}{\end{enumerate}} \nc{\bnum}{\begin{enumerate}[{\rm(i)}]}
\nc{\cl}{\colon} \nc{\seteq}{\mathbin{:=}} \nc{\re}{\mathrm{re}}
\nc{\im}{\mathrm{im}} \nc{\ran}{\rangle} \nc{\lan}{\langle}
\nc{\on}{\operatorname}

\nc{\Hom}{\on{Hom}} \nc{\Oint}{\mathcal{O}_{\mathrm{int}}}
\nc{\Wt}{\on{Wt}} \nc{\pP}{\widetilde{P}} \nc{\eq}{\begin{eqnarray}}
\nc{\eneq}{\end{eqnarray}} \nc{\eqn}{\begin{eqnarray*}}
\nc{\eneqn}{\end{eqnarray*}} \nc{\Lemma}{\begin{lem}}
\nc{\enlemma}{\end{lem}}

\nc{\hs}{\hspace*} \nc{\bfi}{\mathbf{i}} \nc{\eps}{\varepsilon}
\nc{\ba}{\begin{array}} \nc{\ea}{\end{array}}

\nc{\tf}{\tilde{f}} \nc{\id}{\operatorname{id}} \nc{\bl}{\bigl}
\nc{\br}{\bigr} 

\begin{document}

\title[Polyhedral realization of Crystal bases for GKM Algebras]
 {Polyhedral realization of Crystal bases for Generalized Kac-Moody Algebras}

\author[Dong-Uy Shin]
{Dong-Uy Shin}

\address{Department of Mathematics\\
         Chonnam National University\\
         Kwangju 500-757, Korea}
\email{dushin@chonnam.ac.kr}


\subjclass[2000]{Primary 81R50, Secondary 17B37}

\keywords{crystal basis, generalized Kac-Moody algebras, Monster
algebras}

\maketitle

\begin{abstract}
In this paper, we give polyhedral realization of the crystal
$B(\infty)$ of $U_q^-(\mathfrak g)$ for the generalized Kac-Moody
algebras. As applications, we give explicit descriptions of crystals
for the generalized Kac-Moody algebras of rank 2, 3 and {\it Monster
Lie algebras}.
\end{abstract}


\section*{Introduction}
\label{intro}

In his study of Conway and Norton's {\it Moonshine Conjecture}
\cite{CN} for the infinite dimensional ${\bf Z}$-graded
representation $V^\natural$ of the Monster sporadic simple group,
Borcherds introduced a new class of infinite dimensional Lie
algebras called the {\it generalized Kac-Moody algebras}
\cite{Bor1,Bor2}. The structure and representation theories of
generalized Kac-Moody algebras are very similar to those of
Kac-Moody algebras, and a lot of facts about Kac-Moody algebras can
be extended to generalized Kac-Moody algebras. The main difference
is that the generalized Kac-Moody algebras may have simple roots
with non-positive norms whose multiplicity can be greater than one,
called {\it imaginary simple roots}, and they may have infinitely
many simple roots.

The {\it quantum groups} $U_q(\frak g)$ introduced by Drinfel'd
and Jimbo, independently are $q$-deformations of the universal
enveloping algebras $U(\frak g)$ of Kac-Moody algebras $\frak g$
\cite{Drin,Jim}. The important feature of quantum groups is that
the representation theory of $U_q(\frak g)$ is the same as that of
$U(\frak g)$.
Therefore, to understand the structure of representations over
$U_q(\frak g)$, it is enough to understand that of representations
over $U_q(\frak g)$ for some special parameter $q$ which is easy
to treat. The {\it crystal basis theory} which can be viewed as
the representation theory at $q=0$ was introduced by Kashiwara
\cite{Kas91}. Among others, he showed that there exist a crystal
basis $B(\infty)$ for the negative part of a quantum group and a
crystal basis $B(\la)$ for the irreducible highest weight module
$V(\la)$ with a dominant integral highest weight $\la$. Crystal
bases are given a structure of colored oriented graphs, called the
{\it crystal graphs}, which reflect the combinatorial structure of
integrable modules.
So one of the most fundamental problems in the crystal basis theory
is to construct the crystal basis explicitly. In many articles, one
can find several kinds of realizations of crystal bases using
combinatorial objects (for example,
\cite{Kang03,KKM,KMN1,KKS1,KasNak,Li2,Nakashima,NZ}).

In \cite{Kang95}, Kang introduced the {\it quantum generalized
Kac-Moody algebras} $U_q(\frak g)$ -- the quantum groups associated
with generalized Kac-Moody algebras $\frak g$, and he also showed
that, for a generic $q$, the Verma modules and the unitarizable
highest weight modules over $\g$ can be deformed to those over
$U_q(\g)$. 
In \cite{JKK}, Jeong, Kang and Kashiwara developed the crystal basis
theory for quantum generalized Kac-Moody algebras. As in the
Kac-Moody algebra case, they showed that there exist a crystal basis
$B(\infty)$ for the negative part of a quantum generalized Kac-Moody
algebra and a crystal basis $B(\la)$ for the irreducible highest
weight module $V(\la)$ with a dominant integral highest weight
$\la$. However, unfortunately, there is no explicit realization of
crystal bases over quantum generalized Kac-Moody algebras using some
combinatorial objects.

Recently, in \cite{JKKS}, Jeong, Kang, Kashiwara and the author
introduced the notion of {\it abstract crystals} for quantum
generalized Kac-Moody algebras, and the embedding of crystals
$\Psi_{\iota}: B(\infty)\hookrightarrow{\bf Z}_{\ge
0,\iota}^{\infty}$ where $\iota$ is an infinite sequence from the
index set of simple roots. This embedding $\Psi_{\iota}$ of
crystals is  an analogue of the crystal embedding (in Kac-Moody
case) introduced by Kashiwara \cite{Kas93}. But, as in the
Kac-Moody case, in general, it is not easy to find the image ${\rm
Im} \Psi_{\iota}$. In this paper, we give explicit description of
${\rm Im} \Psi_{\iota}$ by a unified method introduced by
Nakashima and Zelevinsky \cite{NZ}, called the {\it polyhedral
realization}. The main obstacle to apply Nakashima and
Zelevinsky's method to the quantum generalized Kac-Moody algebras
was quite different tensor product rule of Kashiwara operators of
crystal bases.

As applications, we give explicit descriptions of the crystals
over Kac-Moody algebras of rank 2 and 3. Finally, for the {\it
Monster Lie algebra} which played an important role in proving the
Moonshine conjecture, we give the explicit description of ${\rm
Im} \Psi_{\iota}$. Since the root multiplicity of Monster Lie
algebra is closely related with the $i$-th coefficient $c(i)$ of
the elliptic modular function $j(q)-744$, we expect that we can
obtain some properties about $c(i)$'s.

\vskip 3mm {\bf Acknowledgments.} \ The author would like to express
his sincere gratitude to Professor S.-J. Kang and Professor J.-A.
Kim for their interest in this work and many valuable discussions.

\section{Crystal Bases for quantum generalized Kac-Moody algebras}
\subsection{Quantum Generalized Kac-Moody Algebras} Let $I$ be a
countable index set. A real matrix $A=(a_{ij})_{i,j\in I}$ is
called a {\it Borcherds-Cartan matrix} if it satisfies: (i)
$a_{ii}=2$ or $a_{ii}\le 0$ for all $i\in I$, (ii) $a_{ij}\le 0$
if $i\neq j$, (iii) $a_{ij}\in{\bf Z}$ if $a_{ii}=2$, (iv)
$a_{ij}=0$ if and only if $a_{ji}=0$. Let $I^{re}=\{i\in I \mid
a_{ii}=2\}$ and $I^{im}=\{i\in I \mid a_{ii}\le 0\}$. Moreover, we
say that an index $i$ in $I^{re}$ (resp. $I^{im}$) is {\it real}
(resp. {\it imaginary}).

In this paper, we assume that for all $i,j\in I$, $a_{ij}\in {\bf
Z}$, $a_{ii}\in 2{\bf Z}$, and $A$ is {\it symmetrizable}. That is,
there is a diagonal matrix $D={\rm diag}(s_i\in {\bf Z}_{>0}|\,i\in
I)$ such that $DA$ is symmetric. We set a Borcherds-Cartan datum
($A,P^{\vee},P, \Pi^{\vee}, \Pi$) as follows:
\begin{equation*}
\begin{array}{l}
A: \,\,\,\text{a Borcherds-Cartan matrix},\\
P^{\vee}=\Big(\bigoplus_{i\in I}{\bf
Z}h_i\Big)\oplus\Big(\bigoplus_{i\in I}{\bf Z}d_i\Big):\,\,\,\text{a free abelian group},\\
P=\{\lambda\in \frak h^*\mid \lambda(P^{\vee})\subset {\bf Z}\}:\,\,\,\text{\it the weight lattice},\\
\Pi^{\vee}=\{h_i\mid i\in I\}\subset \frak h:\,\,\,\text{the set of {\it simple coroots},}\\
\Pi=\{\alpha_i\mid i\in I\}\subset \frak h^*:\,\,\,\text{the set
of
{\it simple roots}.}\\
\end{array}
\end{equation*}
Here, the simple roots $\alpha_i$ ($i\in I$) are defined by
\begin{equation*}
\langle h_j, \alpha_i\rangle=a_{ji}\,\,\,\text{and}\,\,\,\langle
d_j, \alpha_i\rangle=\delta_{ji}.
\end{equation*}
We denote by $P^{+}=\{\la \in P\,|\, \la(h_i) \ge 0 \ \ \text{for
all} \ i\in I \}$ the set of {\it dominant integral weights}. We
also use the notation $Q=\bigoplus_{i\in I} \Z \alpha_i$ and
$Q_{+}=\sum_{i\in I} \Z_{\ge 0} \alpha_i$.

For an indeterminate $q$, set $q_i=q^{s_i}$ and define
\begin{equation*}
[n]_i=\frac{q_i^n-q_i^{-n}}{q_i-q_i^{-1}},\quad [n]_i
!=\prod_{k=1}^{n} [k]_i, \quad {\begin{bmatrix} m \\ n
\end{bmatrix}}_{i} =\frac{[m]_{i}!}{[n]_{i}! \, [m-n]_{i}!}.
\end{equation*}

The {\it quantum generalized Kac-Moody algebra} $U_q({\mathfrak
g})$ associated with a Borcherds-Cartan datum
$(A,P^{\vee},P,\Pi^{\vee},\Pi)$ is the associative algebra over
${\bf Q}(q)$ with $1$ generated by the elements
 $e_i$, $f_i$ $ (i \in I)$ and $q^h$ $(h \in P^{\vee})$ with the
following defining relations:
\begin{equation*}
\begin{aligned}
\ & q^0 =1, \ \ q^h q^{h'} = q^{h+ h'} \quad \text{for\,\, $h, h'\in P^{\vee}$},\\
\ & q^h e_i q^{-h} = q^{\ali(h)} e_i, \quad
q^h f_i q^{-h} = q^{-\ali(h)} f_i \quad \text{for\,\, $h\in P^{\vee}, i\in I$}, \\
\ & e_i f_j - f_j e_i = \delta_{ij} \frac{q^{s_ih_i} -
q^{-s_ih_i}}{q_i - q_i^{-1}} \quad \text{for\,\,\,$i,j \in I$},
\end{aligned}
\end{equation*}

\begin{equation*}
\begin{aligned}
\ & \sum_{k=0}^{1-a_{ij}} (-1)^k {\begin{bmatrix} 1-a_{ij} \\ k
\end{bmatrix}}_{i}
e_i^{1-a_{ij}-k} e_j e_i^{k} = 0 \quad \text{for\,\,\,$a_{ii}=2$, $i \neq j$}, \\
\ & \sum_{k=0}^{1-a_{ij}} (-1)^k {\begin{bmatrix} 1-a_{ij} \\ k
\end{bmatrix}}_{i} f_i^{1-a_{ij}-k} f_j f_i^{k} = 0 \quad \text{for\,\,\,$a_{ii}=2$, $i
\neq j$,}\\
\ &e_ie_j-e_je_i=f_if_j-f_jf_i=0\quad\text{if $a_{ij}=0$.}
\end{aligned}
\end{equation*}
Let us denote by $U_q^+(\g)$ (resp.\ $U_q^-(\g)$) the subalgebra of
$U_q(\g)$ generated by the $e_i$'s (resp.\ the $f_i$'s).


\subsection{Crystal Bases}
The {\it category ${\mathcal O}_{int}$} consists of
$U_q(\g)$-modules $M$ satisfying the following properties:

\begin{itemize}
\item [(i)] $M = \bigoplus_{\mu \in P} M_{\mu}$, where $ M_{\mu} = \{
v\in M \mid q^h v = q^{\mu(h)} v \text{ for all } h\in \pv \} $ is
finite dimensional,

\item [(ii)] there exist finitely many elements $\lambda_1,\cdots,
\lambda_s\in P$ such that $\text{wt}(M)\subset
\bigcup_{j=1}^{s}(\lambda_j-Q_+),$ where $\text{wt}(M)=\{\mu\in P|
M_{\mu}\neq 0\}$,

\item [(iii)] if $a_{ii}=2$, then the action of $f_i$ on $M$ is {\it locally
nilpotent}, i.e., for any $m\in M$ there exists a positive integer
$N$ such that $\fit^N m=0$,

\item [(iv)] if $a_{ii}\le 0$, then $\mu(h_i)\in {\bf Z}_{\ge 0}$ for
every $\mu\in \wt(M)$,

\item [(v)] if $a_{ii}\le 0$ and $\mu(h_i)=0$, then $\fit M_{\mu}=0$,

\item [(vi)] if $a_{ii}\le 0$ and $\mu(h_i)\le -a_{ii}$, then $\eit
M_{\mu}=0$.
\end{itemize}
For instance, the {\it irreducible highest weight module}
$V(\la)=U_q(\g)u_\la$ with $\la\in P^+$ defined by following
relations
\begin{itemize}
\item [(i)] $u_\la$ has weight $\la$,
\item [(ii)] $e_iu_\la=0$ for all $i\in I$,
\item [(iii)] $f_i^{\langle h_i,\la\rangle+1}u_\la=0$ for any $i\in I^{re}$,
\item [(iv)]
$f_iu_\la=0$ if $i\in I^{im}$ and $\langle h_i,\la\rangle=0$,
\end{itemize}
belongs to ${\mathcal O}_{int}$. Moreover, the category ${\mathcal
O}_{int}$ is {\it semisimple} and every simple object in
${\mathcal O}_{int}$ is isomorphic to the irreducible highest
weight module $V(\la)$ with $\la\in P^+$ \cite{JKK}.


Fix an index $i\in I$ and  for $k\ge 0$, set $f_i^{(k)}=f_i^k /
[k]_i!$ if $i$ is real, $f_i^k$ if $i$ is imaginary.
Let $M$ be a $U_q(\frak g)$-module in ${\mathcal O}_{int}$. It was
shown in \cite{JKK} that every weight vector $v\in M_{\lambda}$ can
be written uniquely as
\begin{equation*}
v = \sum_{k\geq0} f_i^{(k)} v_k,
\end{equation*}
where (i) $v_k \in \ker e_i \cap M_{\la+k\ali}$, (ii) if $a_{ii}=2$
and $\langle h_i, \la+n\alpha_i\rangle<n$, then $v_n=0$, and (iii)
if $a_{ii}\le 0$, $n>0$ and $\langle h_i, \la+n\alpha_i\rangle=0$,
then $v_n=0$.
%
%
This expression is called the {\it $i$-string decomposition} of $v$.
The {\it Kashiwara operators} $\eit$ and $\fit$ on $M$ are defined
by
\begin{equation*}
\eit v = \sum_{k\geq1} f_i^{(k-1)} v_k, \qquad \fit v =
\sum_{k\geq0} f_i^{(k+1)} v_k.
\end{equation*}

Let $ \A_0 = \{ f/g \in \Q(q)  \, | \, f, g \in \Q[q], \, g(0)\neq
0 \}$ be the localization of $\Q(q)$ at $(q)$. A {\it crystal
basis} of $M$ is a pair $(L, B)$ such that
\begin{itemize}
\item [(i)] $L$ is a free $\A_0$-submodule of $M$ such that $M \cong
\Q(q) \otimes_{\A_0} L$,

\item [(ii)] $B$ is a $\Q$-basis of $L/qL \cong \Q \otimes_{\A_0} L$,

\item  [(iii)] $L=\bigoplus_{\lambda \in P} L_{\la}$, where $L_{\la} = L
\cap M_{\la}$,

\item  [(iv)] $B=\bigsqcup_{\la \in P} B_{\la}$, where $B_{\la}=B \cap
\left(L_{\la} / q L_{\la} \right)$,

\item  [(v)] $\eit L \subset L$, $\fit L \subset L$ for all $i\in I$,

\item [(vi)] $\eit B\subset B\cup \{0\}$, \ $\fit B\subset B\cup \{0\}$
for all $i\in I$,

\item [(vii)] for all $b, b'\in B$ and $i\in I$, $\fit b = b'$ if and only
if $b = \eit b'$.
\end{itemize}

It was proved in \cite{JKK} that every $M\in {\mathcal O}_{int}$ has
a crystal basis unique up to an automorphism. For $\la\in P^+$,
there is a unique crystal basis $(L(\la),B(\la))$ of $V(\la)$, where
\begin{equation*}
\aligned L(\la)&=\A_0\text{-span}\,\{\tilde f_{i_1} \cdots \tilde
f_{i_r} v_{\la}\,|\,i_k \in I, r\in \Z_{\ge 0}\},\\
B(\la) &= \{ \tilde f_{i_1} \cdots \tilde f_{i_r} v_{\la}+q L(\la)
\in L(\la) / q L(\la) \} \setminus \{0\}. \endaligned
\end{equation*}

\vskip 2mm Fix $i\in I$. For any $P\in U_q^-(\frak g)$, there are
$Q,R\in U_q^-(\frak g)$ such that
$$e_iP-Pe_i=\frac{q^{s_ih_i}Q-q^{-s_ih_i}R}{q_i-q_i^{-1}}.$$
We define the endomorphisms $e_i',e_i'': U_q^-(\frak g)\rightarrow
U_q^-(\frak g)$ by
$$e_i'(P)=R, \quad e_i''(P)=Q.$$
Then every $u\in U_q^-(\frak g)$ can be written uniquely as
\begin{equation*}
u=\sum_{k\ge 0}f_i^{(k)}u_k,
\end{equation*}
where $e_i'u_k=0$ for all $k\ge 0$ and $u_k=0$ for $k\gg 0$.
Moreover, we have $u_k=q_i^{a_{ii}k(k-1)/4}P_i e_i^{(k)}u$, which is
called the {\it $i$-string decomposition} of $u$ \cite{JKK}. The
Kashiwara operators $\eit,\fit$ on $U_q^-(\frak g)$ are defined by
\begin{equation*}
\eit u = \sum_{k\geq1} f_i^{(k-1)} u_k, \qquad \fit u =
\sum_{k\geq0} f_i^{(k+1)} u_k.
\end{equation*}
The crystal basis of $U_q^-(\frak g)$ is a pair $(L,B)$ such that
\begin{itemize}
\item [(i)] $L$ is a free $\A_0$-submodule of $U_q^-(\frak g)$ such
that $U_q^-(\frak g)\cong \Q(q)\otimes_{\A_0}L$,

\item  [(ii)] $B$ is a $\Q$-basis of $L/qL\cong \Q\otimes_{\A_0}L$,

\item  [(iii)] $\eit L\subset L$, $\fit L\subset L$ for all $i\in I$,

\item  [(iv)] $\eit B\subset B\cup \{0\}$, $\fit B\subset B\cup \{0\}$
for all $i\in I$,

\item  [(v)] for all $b, b'\in B$ and $i\in I$, $\fit b = b'$ if and
only if $b = \eit b'$.
\end{itemize}

It was proved in \cite{JKK} that there is a unique crystal basis
$(L(\infty),B(\infty))$ of $U_q^-(\frak g)$, where 
\begin{equation*}
\aligned L(\infty)&=\A_0\text{-span}\,\{\tilde f_{i_1} \cdots \tilde
f_{i_r}\cdot 1\,|\,i_k \in I, r\in \Z_{\ge 0}\},\\
B(\infty) &= \{ \tilde f_{i_1} \cdots \tilde f_{i_r}\cdot 1+q
L(\infty) \in L(\infty) / q L(\infty) \} \setminus \{0\}.
\endaligned
\end{equation*}

\section{Abstract Crystals}
In this section, we recall the notion of abstract crystals and their
examples introduced in \cite{JKKS}. Moreover, we introduce a crystal
${\bf Z}_{\ge 0,\iota}^{\infty}$ associated with an infinite
sequence $\iota$.

\subsection{Abstract Crystals}  An {\it abstract crystal} for $U_q(\g)$ or a {\it
$U_q(\g)$-crystal} is a set $B$ together with the maps $\wt : B
\rightarrow P$, $\eit, \fit: B \rightarrow B \cup \{0\}$ $(i\in I)$,
and $\varepsilon_i, \varphi_i: B\rightarrow \Z\cup\{-\infty\}$
$(i\in I)$ such that for all $b\in B$, we have

\begin{itemize}
\item [(i)] $\wt(\eit b) = \wt b + \alpha_i$ if $i\in I$ and $\eit b \neq 0$,

\item [(ii)] $\wt(\fit b) = \wt b - \alpha_i$ if $i\in I$ and
$\fit b \neq 0$,

\item [(iii)] for any $i \in I$ and $b\in B$, $\vphi_i(b) = \vep_i(b) +
\langle h_i, \wt b \rangle$,

\item [(iv)] for any $i\in I$ and $b,b'\in B$,
$\fit b = b'$ if and only if $b = \eit b'$,\label{cond7}

\item [(v)] for any $i \in I$ and $b\in B$
such that $\eit b \neq 0$, we have
\begin{itemize}
\item [(a)]
$\vep_i(\eit b) = \vep_i(b) - 1$, $\vphi_i(\eit b) = \vphi_i(b) + 1$
if $i\in I^{re}$,
\item [(b)]
$\vep_i(\eit b) = \vep_i(b)$ and $\vphi_i(\eit b) = \vphi_i(b) +
a_{ii}$ if $i\in I^{im}$,
\end{itemize}

\item [(vi)] for any $i \in I$ and $b\in B$ such that $\fit b \neq 0$,
we have
\begin{itemize}
\item [(a)]
$\vep_i(\fit b) = \vep_i(b) + 1$ and $\vphi_i(\fit b) = \vphi_i(b) -
1$ if $i\in I^{re}$,
\item [(b)]
$\vep_i(\fit b) = \vep_i(b)$ and $\vphi_i(\fit b) = \vphi_i(b) -
a_{ii}$ if $i\in I^{im}$,
\end{itemize}

\item [(vii)] for any $i \in I$ and $b\in B$ such that $\vphi_i(b) = -\infty$, we
have $\eit b = \fit b = 0$.
\end{itemize}

%
Let $B_1$ and $B_2$ be
crystals. A {\em morphism of crystals} or a {\em crystal morphism}
$\psi: B_1\rightarrow B_2$ is a map $\psi: B_1\to B_2$ such that

\begin{itemize}
\item [(i)] $\wt(\psi(b))=\wt(b)$ for all $b\in B_1$,

\item [(ii)] $\vep_i(\psi(b))=\vep_i(b)$,
$\varphi_i(\psi(b))=\varphi_i(b)$ for all $b\in B_1$, $i\in I$,

\item [(iii)] if $b\in B_1$ and $i\in I$ satisfy $\fit b\in B_1$, then we have
$\psi(\fit b)=\fit\psi(b)$. \label{cond:crysmor2}
\end{itemize}
For a morphism of crystals $\psi: B_1\rightarrow B_2$, $\psi$ is
called a {\it strict morphism} if
\begin{center}
$\psi(\eit b)=\eit\psi(b)$,\,\, $\psi( \fit b)=\fit\psi(b)$\quad for
all $i\in I$ and $b\in B_1$.
\end{center}
Here we understand $\psi(0)=0$. Moreover, $\psi$ is called an {\em
embedding} if the underlying map $\psi: B_1\rightarrow B_2$ is
injective. In this case, we say that $B_1$ is a {\em subcrystal} of
$B_2$. If $\psi$ is a strict embedding, we say that $B_1$ is a {\em
full subcrystal} of $B_2$.

\begin{example}
{\rm (a)} The crystal basis $B(\la)$ of the irreducible highest
weight module $V(\la)$ is an abstract crystal, where the maps
$\vep_i, \varphi_i$ ($i\in I$) are given by

\begin{equation*}
\text{$\vep_i (b)=\max \{k\ge 0| \eit^k b \neq 0\}$,\,\, $\vphi_i
(b)=\max \{k\ge 0| \fit^k b \neq 0\}$\,\, for $i\in
I^{re}$,}\\
\end{equation*}
\begin{equation*}
\text{$\vep_i
(b)=0$,\,\,$\vphi_i (b)=\langle h_i,\wt(b)\rangle$\,\, for $i\in
I^{im}$.}
\end{equation*}


{\rm (b)} The crystal basis $B(\infty)$ of $U_q^{-}(\g)$ is an
abstract crystal, where the maps $\vep_i, \varphi_i$ ($i\in I$) are
given by
\begin{equation*}
\begin{aligned}
\vep_i (b) & =
\begin{cases}
\max \{k\ge 0| \eit^k b \neq 0\} &\text{for $i\in I^{re}$,} \\
0&\text{for $i\in I^{im}$,}\end{cases}\\
\vphi_i (b) & = \vep_i (b)+\langle h_i,\wt(b)\rangle \quad\text{for
$i\in I$.}
\end{aligned}
\end{equation*}
\end{example}

\begin{example}
For $i\in I$, let $B_{i}=\{b_{i}(-n) \mid n \ge 0 \}$ and define
\begin{equation*}
\begin{aligned}
& \wt (b_i(-n)) = -n \alpha_i, \\
& \eit b_i(-n) = b_i(-n+1), \quad \fit b_i(-n) = b_i(-n-1), \\
& \tilde e_j b_i(-n) = \tilde f_j b_i(-n) = 0 \quad \text{if} \ \
j \neq i,\\
& \vep_i (b_i(-n)) = n, \quad \vphi_i (b_i(-n))=-n \quad \text{if}
\ \ i\in I^{re}, \\
& \vep_i (b_i(-n)) = 0, \quad \vphi_i (b_i(-n))=\langle
h_i,\wt(b_i(-n))\rangle=-na_{ii} \quad \text{if}
\ \ i\in I^{im}, \\
& \vep_j (b_i(-n)) = \vphi_j (b_i(-n)) = -\infty \quad \text{if} \ j
\neq i.
\end{aligned}
\end{equation*}
Here, we understand $b_i(-n)=0$ for $n<0$. Then $B_{i}$ is an
abstract crystal, and it is called an {\it elementary crystal}
\cite{JKKS}.
\end{example}

We define the tensor product of a pair of crystals as follows: for
two crystals $B_1$ and $B_2$, their tensor product $B_1\otimes B_2$
is $\{b_1\otimes b_2\,|\, b_1\in B_1, b_2\in B_2\}$ with the
following crystal structure. The maps $\wt, \vep_i,\varphi_i$ are
given by
\begin{equation*}
\aligned
\wt(b\otimes b')&=\wt(b)+\wt(b'),\\
\vep_i(b\otimes b')&=\max(\vep_i(b), \vep_i(b')-\langle h_i, \wt(b)
\rangle),\\
\varphi_i(b\otimes b')&=\max(\varphi_i(b)+\langle h_i, \wt(b')
\rangle,\varphi_i(b')).
\endaligned
\end{equation*}
For $i\in I$, we define
\begin{equation*}
\fit(b\otimes b')=
\begin{cases}\fit b\otimes b'
&\text{if $\varphi_i(b)>\vep_i(b')$,}\\
b\otimes \fit b' &\text{if $\varphi_i(b)\le \vep_i(b')$.}
\end{cases}
\end{equation*}
For $i\in I^{re}$, we define
\begin{equation*}
\eit(b\otimes b')=
\begin{cases}\eit b\otimes b'\ &\text{if
$\varphi_i(b)\ge \vep_i(b')$,}\\
b\otimes \eit b' &\text{if $\varphi_i(b)< \vep_i(b')$,}
\end{cases}
\end{equation*}
and, for $i\in I^{im}$, we define
\begin{equation*}
\eit(b\otimes b')=
\begin{cases}\eit b\otimes b'\
&\text{if $\varphi_i(b)>\vep_i(b')-a_{ii}$,}\\
0&\text{if $\vep_i(b')<\varphi_i(b)\le\vep_i(b') -a_{ii}$,}\\
b\otimes \eit b' &\text{if $\varphi_i(b)\le\vep_i(b')$.}
\end{cases}
\end{equation*}
This tensor product rule is different from the one given in
\cite{JKK}. But when $B_1=B(\la)$ and $B_2 = B(\mu)$ for $\la, \mu
\in P^{+}$, the two rules coincide. Note that by the definition
above,  $B_1\otimes B_2$ is a crystal. Moreover, it is easy to see
that the associativity law for the tensor product holds \cite{JKKS}.


%
%

\subsection{Crystal structure of ${\bf Z}_{\ge 0,\iota}^{\infty}$}
Let $\iota=(\dots,i_k,\dots,i_1)$ be an infinite sequence such that
\begin{equation}\label{eq:seq-con}
i_k\neq i_{k+1}\quad\text{and}\quad \#\{k\,|\,
i_k=i\}=\infty\quad\text{for any $i\in I$.}
\end{equation}
Now, we give a crystal structure ${\bf Z}_{\ge 0,\iota}^{\infty}$ on
the set of infinite sequences of nonnegative integers
\begin{equation*}
{\bf Z}_{\ge 0}^{\infty}:=\{(\dots,x_k,\dots,x_1)\,|\, x_k\in {\bf
Z}_{\ge 0}\,\,\,\text{and $x_k=0$ for $k\gg 0$}\}
\end{equation*}
associated with $\iota$ as follows: Let
$\overrightarrow{x}=(\dots,x_k,\dots,x_1)$ be an element of ${\bf
Z}_{\ge 0}^{\infty}$. For $k\ge 1$, we define
\begin{equation}\label{eq:sigmak}
\sigma_k(\overrightarrow{x})=\begin{cases} x_k+\sum_{j>k}\langle
h_{i_k},\alpha_{i_j}\rangle x_j &\text{if $i_k\in
I^{re}$,}\\
\sum_{j>k}\langle h_{i_k},\alpha_{i_j}\rangle x_j&\text{if $i_k\in
I^{im}.$}
\end{cases}
\end{equation}
Let
\begin{equation*}
\aligned &\sigma^{(i)}(\overrightarrow{x})=\text{max}_{k:
i_k=i}\{\sigma_k(\overrightarrow{x})\},\\
&n_f=\text{min}\{k\,|\, i_k=i,\,\,
\sigma_k(\overrightarrow{x})=\sigma^{(i)}(\overrightarrow{x})\},\\
&n_e=\begin{cases}\text{max}\{k\,|\, i_k=i,\,\,
\sigma_k(\overrightarrow{x})=\sigma^{(i)}(\overrightarrow{x})\}&\text{if\, $i\in I^{re}$,}\\
n_f &\text{if\, $i\in I^{im}$.}
\end{cases}
\endaligned
\end{equation*}

Now, we define
\begin{equation*}
\fit \overrightarrow{x}=(x_k+\delta_{k,n_f})_{k\ge 1},
\end{equation*}
and
\begin{equation}\label{eq:eit}
\eit \overrightarrow{x}=\begin{cases}(x_k-\delta_{k,n_e})_{k\ge
1}&\text{if $\overrightarrow{x}$ satisfies the condition {\bf
(EC)}},\\
0&\text{otherwise,}
\end{cases}
\end{equation}
where, the condition {\bf (EC)} is as follows:
\begin{equation*}
\aligned {\bf (EC)}\quad &{\rm(i)}\,i\in I^{re}:
\sigma^{(i)}(\overrightarrow{x})>0,\\
&{\rm(ii)}\, i\in I^{im}:
\text{for $k=n_e$ with $k^{(-)}\neq 0$,}\\
&\qquad\qquad\text{ $x_k>1$,\, or $x_k=1$ and
$\sum_{k^{(-)}<j<k}\langle h_i, \alpha_{i_j}\rangle x_j<0$.}
\endaligned
\end{equation*}
Here, $k^{(-)}$ is the maximal index $j<k$ such that $i_j=i_k$. We
also define
\begin{equation*}
\wt(\overrightarrow{x})=-\sum_{j=1}^{\infty}x_j\alpha_{i_j},\quad
\varepsilon_i(\overrightarrow{x})=\sigma^{(i)}(\overrightarrow{x}),\quad
\varphi_i(\overrightarrow{x})=\langle
h_i,\wt(\overrightarrow{x})\rangle+\varepsilon_i(\overrightarrow{x}).
\end{equation*}
 It is easy to see that ${\bf
Z}_{\ge 0}^{\infty}$ is a crystal. We denote this crystal by ${\bf
Z}_{\ge 0,\iota}^{\infty}$.

\begin{remark}
Since $x_k=0$ for $k\gg 0$, it is clear that
$\varepsilon_i(\overrightarrow{x})=0$ for each $i\in I^{im}$, and so
$\varphi_i(\overrightarrow{x})=\langle
h_i,\wt(\overrightarrow{x})\rangle$.
\end{remark}

\subsection{Embedding of Crystals}

\begin{prop}\label{prop:embedding} {\rm \cite{JKKS}}
For all $i \in I$, there exists a unique strict embedding
\begin{equation*}
\Psi_i : B(\infty) \longrightarrow B(\infty) \otimes B_i \quad
\text{such that} \quad u_{\infty} \mapsto u_{\infty} \otimes b_i(0),
\end{equation*}
where $u_{\infty}$ is the highest weight vector in $B(\infty)$.
\end{prop}

%
%

The Proposition \ref{prop:embedding} yields a procedure to determine
the structure of the crystal $B(\infty)$ in terms of elementary
crystals. Take an infinite sequence ${\iota} =(\dots, i_2,i_1)$ in
$I$ such that every $i\in I$ appears infinitely many times. 
For each $N\ge 1$,
taking the composition of crystal embeddings repeatedly, we obtain a
strict crystal embedding
\begin{equation}
\aligned &\Psi^{(N)}\seteq(\Psi_{i_N}\otimes \id\otimes
\cdots\otimes
\id)\circ\cdots\circ(\Psi_{i_2}\otimes \id)\circ \Psi_{i_1}\cl \\
&B(\infty)\hookrightarrow B(\infty)\otimes B_{i_1}\hookrightarrow
B(\infty)\otimes B_{i_2}\otimes B_{i_1}\hookrightarrow\\
&\qquad\qquad\qquad\qquad\qquad\qquad \cdots \hookrightarrow
B(\infty)\otimes B_{i_N}\otimes \cdots\otimes B_{i_1}.
\endaligned
\end{equation}
It is easily seen that, for any $b\in B$, there exists $N>0$ such
that
\begin{equation*}
\Psi^{(N)}(b)=u_{\infty}\otimes b_{i_N}(-x_N)\otimes \cdots\otimes
b_{i_1}(-x_1)
\end{equation*}
for some $x_1,\dots,x_N\in \Z_{\ge 0}$ and $x_k=0$ for $k>N$. Thus
the sequence $(\dots,0,x_N,\dots,x_1)$ belongs to ${\bf Z}_{\ge
0,\iota}^{\infty}$, and so we obtain a map
$$\begin{array}{lclc}
\Psi_{\iota}:&B(\infty)&\to &{\bf Z}_{\ge 0,\iota}^{\infty}\\
&b&\mapsto &(\dots,0,x_N,\dots,x_1).
\end{array}
$$
We can easily see that it is a strict embedding (See also
\cite{JKKS}).

\section{Polyhedral realizations of $B(\infty)$}
In \cite{NZ}, Nakashima and Zelevinsky gave a polyhedral realization
of the crystal base $B(\infty)$ of the negative part $U_q^-(\frak
g)$ of the quantum group $U_q(\frak g)$ associated with Kac-Moody
algebra. In this section, we extend their theory to the case of
quantum generalized Kac-Moody algebras.

\subsection{Polyhedral realizations of $B(\infty)$}

Let $\iota=(i_k)_{k\ge 1}$ be a sequence of indices satisfying
\eqref{eq:seq-con}. Let ${\bf Q}^{\infty}$  be an infinite
dimensional vector space
\begin{center}
${\bf Q}^{\infty}=\{\overrightarrow{x}=(\dots,x_k,\dots,x_1)\mid
x_k\in {\bf Q}\,\,\text{and $x_k=0$ for $k\gg 0$}\}$.
\end{center}
For a linear functional $\psi\in ({\bf Q}^{\infty})^*$, we write
$\psi(\overrightarrow{x})=\sum_{k\ge 1}\psi_k x_k$ ($\psi_k\in {\bf
Q}$). For each $k\ge 1$, we denote by $k^{(+)}$ (resp. $k^{(-)}$)
the minimal (resp. maximal) index $j>k$ (resp. $j<k$) such that
$i_j=i_k$. Let $\beta_k\in ({\bf Q}^{\infty})^*$ be a linear form
\begin{equation}\label{eq:beta}
\aligned
\beta_k(\overrightarrow{x})&=\sigma_k(\overrightarrow{x})-\sigma_{k^{(+)}}(\overrightarrow{x})\\
&=\begin{cases} x_k+\sum_{k<j<k^{(+)}}\langle
h_{i_k},\alpha_{i_j}\rangle x_j+x_{k^{(+)}}&\text{if $i_k\in
I^{re}$,}\\
\sum_{k<j\le k^{(+)}}\langle h_{i_k},\alpha_{i_j}\rangle x_j
&\text{if $i_k\in I^{im}$,}
\end{cases}
\endaligned
\end{equation}
and we set $\beta_0(\overrightarrow{x})=0$. Then, we define a
piecewise-linear operator $S_k=S_{k,\iota}$ on $({\bf
Q}^{\infty})^*$ by
\begin{equation*}
S_k(\psi)=\begin{cases} \psi-\psi_k\beta_k &\text{if $\psi_k>0$,
$i_k\in I^{re}$,}\\
\psi-\psi_k(x_k+\sum_{k<j<k^{(+)}}\langle h_{i_k},
\alpha_{i_j}\rangle x_j-x_{k^{(+)}}) &\text{if $\psi_k>0$, $i_k\in
I^{im}$,}\\
\psi-\psi_k\beta_{k^{(-)}} &\text{if $\psi_k\le 0$.}
\end{cases}
\end{equation*}

Let
\begin{equation*}
\Theta_\iota=\{S_{j_l}\cdots S_{j_1}x_{j_0}\mid l\ge 0,
j_0,\dots,j_l\ge 1\}
\end{equation*} be the set of linear forms obtained from the
coordinate forms $x_j$ by applying transformations $S_k$. Moreover,
for a given $s,t\ge 1$ ($t>s$), let $\Theta^{s\backslash t}_\iota$
be the subset of $\Theta_\iota$ of linear forms obtained from the
coordinate forms $x_s$ by applying transformations $S_k$ with $k\neq
t$, i.e.,
\begin{equation*}
\Theta^{s\backslash t}_\iota=\{S_{j_l}\cdots S_{j_1}x_s\mid l\ge 0,
s,j_1,\dots,j_l\ge 1\},
\end{equation*}
where none of $j_1,\dots,j_l$ is $t$. We impose on $\iota$ the
positivity assumption given in \cite{NZ}. That is,
\begin{equation} \label{eq:posi-assumption}
\text{if $k^{(-)}=0$, \,\,\,then\,\, $\psi_k\ge 0$\,\,for any
$\psi=\sum \psi_jx_j\in \Theta_\iota$.}
\end{equation}


Then we have the following main theorem.

\begin{thm} \label{thm:Main}
Let $\iota$ be a sequence of indices satisfying \eqref{eq:seq-con}
and \eqref{eq:posi-assumption}. Let
$\Psi_{\iota}:B(\infty)\hookrightarrow {\bf Z}_{\ge
0,\iota}^{\infty}$ be the crystal embedding. Then ${\rm
Im}\Psi_\iota$ is the set $\Gamma_{\iota}$ of
$\overrightarrow{x}\in{\bf Z}_{\ge 0,\iota}^{\infty}$ satisfying the
following conditions{\rm :}
\begin{itemize}
\item [{\rm(i)}] $\psi(\overrightarrow{x})\ge 0$ for any $\psi\in
\Theta_\iota$,

\item [{\rm(ii)}] for each $t$ with $i_t\in I^{im}$, if $x_t\neq
0$ and $t^{(-)}\neq 0$, then
\begin{equation}\label{eq:mainii-1}
\sum_{t^{(-)}<j<t}\langle h_{i_t}, \alpha_{i_j}\rangle x_j< 0.
\end{equation}
In addition, if $\langle h_{i_t}, \alpha_{i_j}\rangle x_j=0$
$(t^{(-)}<j<t)$ for all $j$ with $i_j\in I^{im}$, there exists an
integer $p$ $(t^{(-)}<p<t)$ such that $i_p\in I^{re}$,
\begin{equation}\label{eq:mainii-2}
\text{$\langle h_{i_t}, \alpha_{i_p}\rangle x_p<0$ and\,\,
$\psi(\overrightarrow{x})>0$ for any $\psi\in \Theta^{p\backslash
t}_\iota$.}
\end{equation}
\end{itemize}
\end{thm}

\begin{proof}
We prove it in Subsection 3.2 later.
\end{proof}

\begin{cor} \label{cor:allim} Assume that all elements of $I$ are imaginary,
that is, $I=I^{im}$. Then the image of the crystal embedding ${\rm
Im}\Psi_\iota$ equals to the set of $\overrightarrow{x}\in {\bf
Z}_{\ge 0,\iota}^{\infty}$ satisfying \eqref{eq:mainii-1} of Theorem
\ref{thm:Main}.
\end{cor}

\begin{proof} By the simple calculation, it is easy to see
that the set $\Theta_{\iota}$ consists of the linear combinations of
the coordinate forms $x_j$ with nonnegative coefficients, which
completes the proof.
\end{proof}


Now, we consider the case that the cardinality of $I^{re}$ is $1$.
Then it is easy to see that $S_jx_j$ is a linear combinations of
$x_k$'s with nonnegative coefficients except for $i_j\in I^{re}$. If
$i_j\in I^{re}$, then
\begin{center}
$S_jx_j=-\sum_{j<t<j^{(+)}}\langle h_j, \alpha_{i_t}\rangle
x_t-x_{j^{(+)}}$,
\end{center}
and
\begin{itemize}
\item [(i)] if $k=j^{(+)}$, then $S_kS_jx_j$ is $x_j$,
\item [(ii)] if $j<k<j^{(+)}$ and $\langle h_{i_k},\alpha_{i_j}\rangle
<0$, then $S_kS_jx_j$ is a linear combination of $x_t$'s of
nonnegative coefficients,
\item [(iii)] if $k$ does not belong to the cases (i) and (ii), then $S_kS_jx_j$ is $S_jx_j$ itself.
\end{itemize}
Therefore, it is easy to see that the condition (i) of Theorem
\ref{thm:Main} is changed to
\begin{equation}
\text{$S_jx_j\ge 0$\,\, for all $j$ with $i_j\in I^{re}$.}
\end{equation}
Moreover, for given $p,t$ in (ii) of Theorem \ref{thm:Main}, since
$\psi_t>0$ for any $\psi\in \Theta^{p\backslash t}_{\iota}$, the
above (i)-(iii) implies that the condition $S_px_p>0$ is the same as
the condition that $\psi(\overrightarrow{x})>0$ for all $\psi\in
\Theta^{p\backslash t}_{\iota}$. Finally, by the above (i)-(iii), it
is clear that any sequence $\iota$ satisfies the positivity
assumption \eqref{eq:posi-assumption}. Therefore, we have the
following simple and important corollary.

%

\begin{cor}\label{cor:re1}
Let $I$ be an index set such that the cardinality of $I^{re}$ is
$1$, and let $\iota$ be a sequence of indices in $I$ satisfying
\eqref{eq:seq-con}. Then the image ${\rm Im}\Psi_\iota$ of the
crystal embedding is the set $\Gamma_{\iota}$ of
$\overrightarrow{x}\in{\bf Z}_{\ge 0,\iota}^{\infty}$ satisfying the
following conditions:
\begin{itemize}
\item [{\rm(i)}] $S_jx_j\ge 0$ for all $j$ with $i_j\in I^{re}$,

\item [{\rm(ii)}] for each $t$ with $i_t\in I^{im}$, if $x_t\neq
0$ and $t^{(-)}\neq 0$, then
\begin{equation*}
\sum_{t^{(-)}<j<t}\langle h_{i_t}, \alpha_{i_j}\rangle x_j< 0.
\end{equation*}
In addition, if $\langle h_{i_t}, \alpha_{i_j}\rangle x_j=0$
$(t^{(-)}<j<t)$ for all $i_j\in I^{im}$, there exists an integer $p$
$(t^{(-)}<p<t)$ such that $i_p\in I^{re}$,
\begin{equation*}
\text{$\langle h_{i_t}, \alpha_{i_p}\rangle x_p<0$ and\,\,
$S_px_p>0$.}
\end{equation*}
\end{itemize}

%
%

\end{cor}

\begin{example}
Assume that $I=\{1,2\}$ and $\iota=(\dots,2,1,2,1)$. Set
\begin{equation*}
\text{$\alpha_1(h_1)=-a$,\,\, $\alpha_1(h_2)=-c$,\,\,
$\alpha_2(h_1)=-b$\,\, and\,\,$\alpha_2(h_2)=2$}
\end{equation*}
where $a,b,c\in {\bf Z}_{\ge 0}$. Then $I^{re}=\{2\}$,
$I^{im}=\{1\}$, and if $k\ge 3$, $k^{(-)}\neq 0$. Therefore, for
each $k\ge 1$, if $x_{2k+1}\neq 0$, we have
\begin{equation*}
\sum_{2k-1<j<2k+1}\langle h_{1}, \alpha_{i_j}\rangle x_j=-bx_{2k}<
0.
\end{equation*}
Moreover, since $x_{2k}>0$ and $i_{2k}=2\in I^{re}$, we have
$S_{2k}x_{2k}=x_{2k}-\beta_{2k}=cx_{2k+1}-x_{2k+2}>0$. Therefore, by
Corollary \ref{cor:re1} the image of the crystal embedding ${\rm
Im}\Psi_\iota$ is given by the subset $\Gamma_{\iota}$ of
$\overrightarrow{x}\in{\bf Z}_{\ge 0,\iota}^{\infty}$ as follows:

{\rm (a)} When $b=c=0$,
\begin{center}
$x_{k}=0$ for $k\ge 3$.
\end{center}

{\rm (b)} When neither $b$ nor $c$ is $0$,
\begin{equation*}
\aligned
&\text{{\rm(i)} for each $k\ge 1$, $cx_{2k+1}-x_{2k+2}> 0$ unless $x_{2k+1}=x_{2k+2}=0$,} \\
&\text{{\rm(ii)} for each $k\ge 1$, if $x_{2k+1}\neq 0$, then
$x_{2k}>0.$}
\endaligned
\end{equation*}
\end{example}

\subsection{The proof of Theorem \ref{thm:Main}}

We know that $\text{Im}\Psi_\iota$ is a subcrystal of ${\bf Z}_{\ge
0, \iota}^{\infty}$ obtained by applying the Kashiwara operators
$\fit$ to $\Psi_\iota(u_{\infty})=\overrightarrow{0}=(\dots,0,0,0)$
and $\overrightarrow{0}$ belongs to $\Gamma_\iota$. So, in order to
prove that ${\rm Im}\Psi_{\iota}\subset \Gamma_{\iota}$, it suffices
to show that $\Gamma_{\iota}$ is closed under all $\fit$. Let
$\overrightarrow{x}\in \Gamma_{\iota}$ and $i\in I$. Suppose that
$\fit \overrightarrow{x}=(\dots,x_k+1,\dots,x_1)$. Since
\begin{equation*}
\psi(\fit \overrightarrow{x})=\psi(\overrightarrow{x})+\psi_k\ge
\psi_k \quad\text{for any $\psi\in \Theta_{\iota}$,}
\end{equation*}
in order to prove (i), it is enough to consider the case when
$\psi_k<0$. By the positivity condition \eqref{eq:posi-assumption}
of $\iota$, we have $k^{(-)}\ge 1$. By \eqref{eq:sigmak}, we have
$\sigma_k(\overrightarrow{x})>\sigma_{k^{(-)}}(\overrightarrow{x})$
(Indeed, when $i_k\in I^{im}$, $\sigma_k(\overrightarrow{x})=0$ and
$\sigma_{k^{(-)}}(\overrightarrow{x})<0$), and so
\begin{equation*}
\beta_{k^{(-)}}(\overrightarrow{x})=
\sigma_{k^{(-)}}(\overrightarrow{x})-\sigma_k(\overrightarrow{x})\le
-1.
\end{equation*}
Therefore,
\begin{equation}\label{eq:proof-fit} \aligned
\psi(\fit\overrightarrow{x})&=\psi(\overrightarrow{x})+\psi_k\\
&\ge
\psi(\overrightarrow{x})-\psi_k\beta_{k^{(-)}}(\overrightarrow{x})\\
&=(S_k\psi)(\overrightarrow{x})\ge 0.
\endaligned
\end{equation}
Now, suppose that $\fit \overrightarrow{x}$ does not satisfy the
condition \eqref{eq:mainii-1}. Then $k=t$, and
\begin{equation*}
\text{$x_t=0$,\quad $\sum_{t^{(-)}<j<t}\langle h_{i_t},
\alpha_{i_j}\rangle x_j=0$\quad in $\overrightarrow{x}$.}
\end{equation*}
But, it can not occur by the definition of Kashiwara operator
$\fit$. Now, we show that $\fit \overrightarrow{x}$ satisfies the
condition \eqref{eq:mainii-2}. First, suppose that there exist $p$
and $t$ satisfying \eqref{eq:mainii-2} in $\overrightarrow{x}$.
Since
$\psi(\fit\overrightarrow{x})=\psi(\overrightarrow{x})+\psi_k$, it
is enough to consider the cases that $\psi_k<0$. Note that by
definition of the set $\Theta^{p\backslash t}_{\iota}$, $\psi_t>0$
for all $\psi\in \Theta^{p\backslash t}_{\iota}$. So it suffices to
the case that $k\neq t$. If $k\neq t$, then $S_k\psi \in
\Theta^{p\backslash t}_{\iota}$ and so
$\psi(\fit\overrightarrow{x})=\psi(\overrightarrow{x})+\psi_k\ge
(S_k\psi)(\overrightarrow{x})>0.$

%
%

%
Second, suppose that $k=t$, $x_t=0$, and for any $j$ such that
$t^{(-)}<j<t$, $i_j\in I^{re}$, $\langle h_{i_t},\alpha_{i_j}\rangle
x_j<0$, there is a $\psi\in\Theta_{\iota}^{j\backslash t}$ such that
$\psi(\overrightarrow{x})=0$ in $\overrightarrow{x}$. Note that
since $j$ is the index such that $\langle h_{i_t},
\alpha_{i_j}\rangle<0$, we have $\psi_t>0$ for all $\psi\in
\Theta_{\iota}^{j\backslash t}$. Therefore, $\psi(\fit
\overrightarrow{x})=\psi(\overrightarrow{x})+\psi_t\ge \psi_t>0$ for
all $\psi\in \Theta_{\iota}^{j\backslash t}$. Therefore, ${\rm
Im}\Psi_{\iota}\subset \Gamma_{\iota}$.


For the proof of the reverse inclusion $\Gamma_{\iota}\subset {\rm
Im}\Psi_i$, note that for any $\overrightarrow{x}\in {\bf Z}_{\ge
0,\iota}^{\infty}\backslash \{\overrightarrow{0}\}$ satisfying the
condition (ii), there is an $i\in I$ such that $\eit
\overrightarrow{x}\neq 0$. Indeed, for  the largest number $k$ such
that $x_k>0$ in $\overrightarrow{x}$, if $i_k\in I^{re}$, then
$\sigma_k(\overrightarrow{x})=x_k>0$ and so
$\sigma^{(i_k)}(\overrightarrow{x})\ge
\sigma_k(\overrightarrow{x})>0$, which implies
$\tilde{e}_{i_k}\overrightarrow{x}\neq 0$. If $i_k\in I^{im}$, then
$n_f=n_e=k$ by the condition \eqref{eq:mainii-1}, and so we have
$\tilde{e}_{i_k}\overrightarrow{x}\neq 0$.

Since $\Gamma_\iota\subset {\bf Z}_{\ge 0,\iota}^{\infty}$, if
$\Gamma_\iota$ is closed under the Kashiwara operators $\eit$ for
all $i\in I$, then for any $\overrightarrow{x}\in \Gamma_\iota$,
there are $i_1,\dots,i_t\in I$ such that
\begin{equation*}
\tilde{e}_{i_t}\dots
\tilde{e}_{i_1}\overrightarrow{x}=\overrightarrow{0}.
\end{equation*}
Moreover, it means that
\begin{equation*}\tilde{f}_{i_1}\dots
\tilde{f}_{i_t}\overrightarrow{0}=\overrightarrow{x},
\end{equation*}
which implies that $\Gamma_{\iota}\subset {\rm Im}\Psi_i$. So it is
enough to show that $\eit \Gamma_\iota\subset \Gamma_\iota\cup\{0\}$
for all $i\in I$. Let $\overrightarrow{x}\in \Gamma_{\iota}$ and
$i\in I$. Suppose that $\eit
\overrightarrow{x}=(\dots,x_k-1,\dots,x_1)$. Since
\begin{equation*}
\psi(\eit \overrightarrow{x})=\psi(\overrightarrow{x})-\psi_k\ge
-\psi_k \quad\text{for any $\psi\in \Theta_{\iota}$,}
\end{equation*} to prove (i) it suffices to consider the case
when $\psi_k>0$. By \eqref{eq:eit}, we have
\begin{equation*}
\beta_k(\overrightarrow{x})=
\sigma_k(\overrightarrow{x})-\sigma_{k^{(+)}}(\overrightarrow{x})\ge
1 \quad(i\in I^{re})
\end{equation*}
and
\begin{equation*}
x_k+\sum_{k<j<k^{(+)}}\langle h_{i}, \alpha_{i_j}\rangle
x_j-x_{k^{(+)}}\ge 1\quad(i\in I^{im}).
\end{equation*}
Therefore,
\begin{equation}\label{eq:proof-eit}
\aligned \psi(\eit
\overrightarrow{x})&=\psi(\overrightarrow{x})-\psi_k\\
&\ge \begin{cases}
\psi(\overrightarrow{x})-\psi_k\beta_k(\overrightarrow{x})&\text{if
$i\in I^{re}$,}\\
\psi(\overrightarrow{x})-\psi_k(x_k+\displaystyle{\sum_{k<j<k^{(+)}}}\langle
h_{i}, \alpha_{i_j}\rangle x_j-x_{k^{(+)}})&\text{if $i\in I^{im}$,}
\end{cases}\\
&=(S_k\psi)(\overrightarrow{x})\ge 0.
\endaligned
\end{equation}

Now, suppose that $\eit \overrightarrow{x}$ does not satisfy the
condition (ii). First, suppose that $\eit \overrightarrow{x}$ does
not satisfy \eqref{eq:mainii-1}. If $i_k=i\in I^{im}$ and
$t^{(-)}<k<t$, then by the definition of Kashiwara operator $\eit$,
we have $\langle h_i,\alpha_{i_t}\rangle=0$, and so $\langle
h_{i_t},\alpha_{i}\rangle=0$. But, in this case, it is clear that
\eqref{eq:mainii-1} holds in $\eit \overrightarrow{x}$. Second,
suppose that $k$ is the unique index such that $t^{(-)}<k<t$ with
$i_k\in I^{re}$ and $\langle h_{i_t}, \alpha_{i_k}\rangle x_k<0$ in
$\overrightarrow{x}$. In this case, $S_k
x_k(\overrightarrow{x})=x_k-\beta_k>0$ by \eqref{eq:mainii-2}, and
by the definition of Kashiwara operator $\eit$, we have $\beta_k>0$.
Hence, $x_k>\beta_k>0$ and so $x_k>1$. Therefore, $\eit
\overrightarrow{x}$ satisfies \eqref{eq:mainii-1}. So it suffices to
consider the case that $\eit \overrightarrow{x}$ does not satisfy
the condition \eqref{eq:mainii-2}. First, suppose that $k=p$ and
$x_p=1$ in $\overrightarrow{x}$. But, since $S_p
x_p(\overrightarrow{x})=x_p-\beta_p>0$, we have $\beta_p\le 0$. It
contradicts the definition of Kashiwara operator $\eit$.

Second, suppose that $k\neq p$. If $k\neq t$, then by the same
argument in \eqref{eq:proof-eit}, we have
$\psi(\eit\overrightarrow{x})>0$ for all $\psi\in
\Theta^{p\backslash t}_\iota$. So, it suffices to consider the case
that $k=t$ and $x_t>1$. However, in this case,
\begin{center}
$x_t+\sum_{t<j<t^{(+)}}\langle h_i,\alpha_{i_j}\rangle
x_j-x_{t^{(+)}}>1$
\end{center} and so
\begin{equation}
\aligned \psi(\eit
\overrightarrow{x})&=\psi(\overrightarrow{x})-\psi_t\\
&>\psi(\overrightarrow{x})-\psi_t(x_t+\displaystyle{\sum_{t<j<t^{(+)}}}\langle
h_{i}, \alpha_{i_j}\rangle x_j-x_{t^{(+)}})\\
&=(S_t\psi)(\overrightarrow{x})\ge 0.
\endaligned
\end{equation}
Therefore, $\Gamma_{\iota}$ is closed under all $\eit$.

\section{Applications: Rank 3 case and Monster Lie algebra}
In this section, we will give an explicit description of the image
of the Kashiwara embedding for the generalized Kac-Moody algebras of
rank 3 and Monster Lie algebras.

\subsection{Rank 3 case}

Assume that $I=\{1,2,3\}$ and $\iota=(\dots,1,3,2,1)$. 
Consider the case when $1,2\in I^{im}$ and $3\in I^{re}$. Let $A$ be
a Borcherds-Cartan matrix
\begin{equation*}
A=\left(\begin{array}{rrr} -a & -b &-c \\
-d &-e & -f\\
-g &-h &2
\end{array}\right),
\end{equation*}
%
%
%
where $a,b,c,d,e,f,g,h\in {\bf Z}_{\ge 0}$. For each $k\ge 1$, we
have
\begin{equation*}
S_{3k}x_{3k}=gx_{k+1}+hx_{k+2}-x_{k+3},
\end{equation*}
Moreover, since $I^{im}=\{1,2\}$, for each $k$ with $i_k=1,2$,
\begin{equation*}
\sum_{k^{(-)}<j<k}\langle h_{i_k}, \alpha_{i_j}\rangle x_j
=\begin{cases} -bx_{k-2}-cx_{k-1} &\text{if $i_k=1$},\\
-fx_{k-2}-dx_{k-1}& \text{if $i_k=2$.}
\end{cases}
\end{equation*}
Therefore, by Corollary \ref{cor:re1}, we have

%

\begin{cor}
Assume that $1,2\in I^{im}$ and $3\in I^{re}$. The image of the
crystal embedding ${\rm Im}\Psi_\iota$ is given by the subset
$\Gamma_{\iota}$ of $\overrightarrow{x}\in{\bf Z}_{\ge
0,\iota}^{\infty}$ satisfying the following conditions:
\begin{itemize}
\item [{\rm(i)}] $gx_{3k+1}+hx_{3k+2}-x_{3k+3}\ge 0$ for $k\ge 1$,

\item [{\rm(ii)}] for each $k\ge 1$, if $x_{3k+1}>0$ {\rm (resp.} $x_{3k+2}>0${\rm)}, then
\begin{center}
$bx_{3k-1}+cx_{3k}>0$\quad {\rm (resp. $fx_{3k}+dx_{3k+1}>0$)}.
\end{center}
Moreover,  if $bx_{3k-1}=0$ {\rm (resp. $dx_{3k+1}=0$)},
\begin{center}
$gx_{3k+1}+hx_{3k+2}-x_{3k+3}> 0.$
\end{center}
\end{itemize}
\end{cor}

\subsection{Monster Lie algebras} Let $I=\{-1\}\cup {\bf
N}$ and let $A=(-(i+j))_{i,j\in I}$ be a Borcherds-Cartan matrix of
charge $\underline{m}=(c(i)|\,i\in I)$. Here, $c(i)$ is the
coefficient of the elliptic modular function
$$j(q)- 744 = q^{-1} + 196884q + 21493760q^2 + \cdots
= \sum_{i=-1}^{\infty} c(i) q^i. $$ Then we have the associated
generalized Kac-Moody algebra called Monster Lie algebra.

On the other hand, let
\begin{equation*}
\text{$I=\{-1=-1_1\}\cup \{i_t\,|\,i\in {\bf N}, t=1,\dots,c(i)\}$
and $A=(-(i+j))_{p,q\in I}$,}
\end{equation*}
where $p=i_l$ and $q=j_m$ for some $1\le l\le c(i)$ and $1\le m\le
c(j)$. Then the associated generalized Kac-Moody algebra is also
the Monster Lie algebra. From now on, we adopt the latter
exposition of the Monster Lie algebra. Assume that
\begin{equation*}
\begin{aligned} \iota=(\dots,-1,3_{c(3)},\dots,
3_1,&2_{c(2)},\dots,2_1,1_{c(1)},\dots,1_1,-1,\\
&2_{c(2)},\dots,2_1,1_{c(1)},\dots,1_1,-1,1_{c(1)},\dots,1_1,-1).
\end{aligned}
\end{equation*}

Let $I_{(-1)}$ be the set of positive integers $t$ such that
$i_t=-1$, i.e.,
\begin{equation*}
I_{(-1)}=\{1\}\cup\{b(n)=nc(1)+(n-1)c(2)+\dots+c(n)+n+1\,|\,n\in
{\bf N}\},
\end{equation*}
and for any $n\ge 1$, we set
\begin{equation*}
\sigma(n)=c(1)+\dots+c(n).
\end{equation*}

\begin{thm} \label{thm:monster-infty}
The image of the Kashiwara embedding ${\rm Im}\Psi_\iota$ is given
by the subset $\Gamma_{\iota}$ of $\overrightarrow{x}\in{\bf
Z}_{\ge 0,\iota}^{\infty}$ such that
\begin{itemize}
\item [{\rm(i)}] $x_{c(1)+2}=0$, and for each $n\ge 1$,
$$
\sum_{k=1}^n
k(x_{b(n)+\sigma(k)+1}+\cdots+x_{b(n)+\sigma(k+1)})-x_{b(n)+\sigma(n+1)+1}\ge
0,
$$
%
%
%

\item [{\rm(ii)}] for each $k\notin I_{(-1)}$, if $x_k>0$ and
$k^{(-)}\neq 0$, then
$$\sum_{k^{(-)}<j<k}\langle h_{i_k}, \alpha_{i_j}\rangle x_j< 0.$$
Moreover, if $\langle h_{i_k}, \alpha_{i_j}\rangle x_j= 0$ for all
$k^{(-)}<j<k$ with $j\notin I_{-1}$, then there exists $m\ge 1$ such
that $k^{(-)}<b(m)<k$ and
$$
\sum_{k=1}^m
k(x_{b(m)+\sigma(k)+1}+\cdots+x_{b(m)+\sigma(k+1)})-x_{b(m)+\sigma(m+1)+1}>
0,
$$
\end{itemize}
\end{thm}

\begin{proof} 
By simple calculation, we have
\begin{equation*}
S_{1}x_{1}=x_{1}-(x_{1}+\langle h_{-1},\alpha_1\rangle (x_2+\dots+
x_{c(1)+1})+x_{c(1)+2})=-x_{c(1)+2}
\end{equation*}
and for each $n\ge 1$

\begin{equation*}
\aligned
S_{b(n)}x_{b(n)}&=x_{b(n)}-(x_{b(n)}+\sum_{k=1}^{n+1}\langle
h_{-1},\alpha_k\rangle (
x_{b(n)+\sigma(k-1)+1}+\dots+x_{b(n)+\sigma(k)})\\
&\qquad\qquad\qquad\qquad +  x_{b(n)+\sigma(n+1)+1})\\
&=\sum_{k=1}^n
k(x_{b(n)+\sigma(k)+1}+\cdots+x_{b(n)+\sigma(k+1)})-x_{b(n)+\sigma(n+1)+1}.
\endaligned
\end{equation*}
%
%
Moreover, it is also easy to see that $S_jS_kx_k$ for all $j$ is a
linear combination of $x_j$'s with nonnegative coefficients.
Therefore, we have the results.
\end{proof}

Finally, by Theorem \ref{thm:monster-infty}, we have the following
character formula of the negative part $U_q^{-}(\frak g)$ of the
quantum Monster Lie algebra $U_q(\frak g)$.

\begin{cor}
\begin{equation*}
{\rm ch}\,U_q^{-}(\frak g)=\sum_{\overrightarrow{x}\in
\Gamma_{\iota}}e^{\wt(\overrightarrow{x})}=\sum_{\overrightarrow{x}\in
\Gamma_{\iota}}e^{-\sum_{j=1}^{\infty}x_j\alpha_{i_j}}.
\end{equation*}

\end{cor}

\vskip 5mm


\label{}




\end{document}

\oneappendix 
\section{About the bibliography}
References in the bibliography should be listed alphabetically by
the authors' surname(s) and, for the same set of authors, by
publication year. Detailed formatting (italic, etc.) should be
avoided; please concentrate on giving full and clear information.

\affiliationone{
   F. Irst and Second Author\\
   Postal Address should be
      added here, including\\
   Country
   \email{first@university.ac.uk\\
   sauthor@university.ac.uk}}
\affiliationtwo{
   T. Hird\\
   Previous postal address where
     the research was performed and\\
   Country
   \email{hird@university.ac.uk}}
\affiliationthree{~} 
\affiliationfour{%
   Current address:\\
   Present long-term address\\
   Country
   \email{t.hird@institution.edu}}
\end{document}